# Binomial upper bounds on generalized moments and tail probabilities of (super)martingales with differences bounded from above

**Iosif Pinelis**[1]

*Michigan Technological University*

**Abstract:** Let $(S_0, S_1, \dots)$ be a supermartingale relative to a nondecreasing sequence of $\sigma$-algebras $H_{\le 0}, H_{\le 1}, \dots$, with $S_0 \le 0$ almost surely (a.s.) and differences $X_i := S_i - S_{i-1}$. Suppose that $X_i \le d$ and $\mathsf{Var}(X_i | H_{\le i-1}) \le \sigma_i^2$ a.s. for every $i = 1, 2, \dots$, where $d > 0$ and $\sigma_i > 0$ are non-random constants. Let $T_n := Z_1 + \dots + Z_n$, where $Z_1, \dots, Z_n$ are i.i.d. r.v.'s each taking on only two values, one of which is $d$, and satisfying the conditions $\mathsf{E} Z_i = 0$ and $\mathsf{Var} Z_i = \sigma^2 := \frac{1}{n}(\sigma_1^2 + \dots + \sigma_n^2)$. Then, based on a comparison inequality between generalized moments of $S_n$ and $T_n$ for a rich class of generalized moment functions, the tail comparison inequality

$$\mathsf{P}(S_n \ge y) \le c\, \mathsf{P}^{\mathsf{Lin,LC}}(T_n \ge y + \tfrac{h}{2}) \quad \forall y \in \mathbb{R}$$

is obtained, where $c := e^2/2 = 3.694\dots$, $h := d + \sigma^2/d$, and the function $y \mapsto \mathsf{P}^{\mathsf{Lin,LC}}(T_n \ge y)$ is the least log-concave majorant of the linear interpolation of the tail function $y \mapsto \mathsf{P}(T_n \ge y)$ over the lattice of all points of the form $nd + kh$ $(k \in \mathbb{Z})$. An explicit formula for $\mathsf{P}^{\mathsf{Lin,LC}}(T_n \ge y + \frac{h}{2})$ is given. Another, similar bound is given under somewhat different conditions. It is shown that these bounds improve significantly upon known bounds.

## 1. Introduction

To begin with, consider normalized Khinchin-Rademacher sums $\varepsilon_1 a_1 + \dots + \varepsilon_n a_n$, where the $\varepsilon_i$'s are i.i.d. Rademacher random variables (r.v.'s), with $\mathsf{P}(\varepsilon_i = \pm 1) = \frac{1}{2}$, and the $a_i$'s are real numbers such that $a_1^2 + \dots + a_n^2 = 1$. Whittle [27] (cf. Haagerup [10]) established the sharp form

$$\mathsf{E}\, f(\varepsilon_1 a_1 + \dots + \varepsilon_n a_n) \le \mathsf{E}\, f\big(\tfrac{1}{\sqrt{n}}(\varepsilon_1 + \dots + \varepsilon_n)\big) \le \mathsf{E}\, f(Z) \tag{1.1}$$

of Khinchin's inequality [16] for the power moment functions $f(x) = |x|^p$ with $p \ge 3$, where $Z \sim N(0,1)$. An exponential version of inequality (1.1), with the moment functions $f(x) = e^{\lambda x}$ for $\lambda \in \mathbb{R}$, follows from a result of Hoeffding [12]. An immediate corollary of that is the exponential inequality

$$\mathsf{P}\ (\varepsilon_1 a_1 + \dots + \varepsilon_n a_n \ge x) \le e^{-x^2/2} \quad \forall x \ge 0. \tag{1.2}$$

(In fact, Hoeffding [12] obtains more general results; cf. Remark 2.2 below.) This upper bound, $e^{-x^2/2}$, invites a comparison with an "ideal" upper bound, of the form $c\,\mathsf{P}(Z \ge x)$ for some absolute constant $c > 0$, which one might expect to have in view of (1.1). Since $\mathsf{P}(Z \ge x) \sim \frac{1}{x\sqrt{2\pi}} e^{-x^2/2}$ as $x \to \infty$, this comparison suggests that a factor of order $\frac{1}{x}$ (for large $x$) is "lost" in (1.2). It turns out that the cause of this loss is that the class of the exponential moment functions is too small. For $\alpha > 0$, consider the following class of functions $f\colon \mathbb{R} \to \mathbb{R}$:

$$\mathcal{F}_+^{(\alpha)} := \Big\{ f\colon f(u) = \int_{-\infty}^{\infty} (u-t)_+^\alpha \, \mu(\mathrm{d}t) \text{ for some Borel measure } \mu \ge 0 \text{ on } \mathbb{R} \text{ and all } u \in \mathbb{R} \Big\},$$

where $x_+ := \max(0, x)$ and $x_+^\alpha := (x_+)^\alpha$. Define $\mathcal{F}_-^{(\alpha)}$ as the class of all functions of the form $u \mapsto f(-u)$, where $f \in \mathcal{F}_+^{(\alpha)}$. Let $\mathcal{F}^{(\alpha)} := \{ f + g\colon f \in \mathcal{F}_+^{(\alpha)}, g \in \mathcal{F}_-^{(\alpha)} \}$.

[1]Department of Mathematical Sciences, Michigan Technological University, Houghton, Michigan 49931, USA, e-mail: ipinelis@math.mtu.edu

**Proposition 1.1.** *For every natural $\alpha$, one has $f \in \mathcal{F}_+^{(\alpha)}$ iff $f$ has finite derivatives $f^{(0)} := f$, $f^{(1)} := f', \dots, f^{(\alpha-1)}$ on $\mathbb{R}$ such that $f^{(\alpha-1)}$ is convex on $\mathbb{R}$ and $f^{(j)}(-\infty+) = 0$ for $j = 0, 1, \dots, \alpha - 1$.*

The proof of this and other statements (whenever a proof is necessary) is deferred to Section 3.

It follows that, for every $t \in \mathbb{R}$, every $\beta \ge \alpha$, and every $\lambda > 0$, the functions $u \mapsto (u-t)_+^\beta$ and $u \mapsto e^{\lambda(u-t)}$ belong to $\mathcal{F}_+^{(\alpha)}$, while the functions $u \mapsto |u-t|^\beta$ and $u \mapsto \cosh \lambda(u-t)$ belong to $\mathcal{F}^{(\alpha)}$.

**Remark 1.2.** Eaton [5] (cf. [7, 21]) obtained inequality (1.1) for a class of moment functions, which essentially coincides with the class $\mathcal{F}^{(3)}$, as seen from [21, Proposition A.1]. Since the class $\mathcal{F}^{(3)}$ is much richer than the class of exponential moment functions, Eaton [6] conjectured (based on asymptotics and numerics) that his inequality (1.1) for $f \in \mathcal{F}^{(3)}$ implies the inequality $\mathsf{P}\,(\varepsilon_1 a_1 + \dots + \varepsilon_n a_n \ge x) \le \frac{2e^3}{9} \frac{1}{x\sqrt{2\pi}} e^{-x^2/2}$ for all $x > \sqrt{2}$, so that the "lost" factor $\frac{1}{x}$ would be restored. A stronger form of this conjecture was proved by Pinelis [21]:

$$\mathsf{P}\,(\varepsilon_1 a_1 + \dots + \varepsilon_n a_n \ge x) \le \frac{2e^3}{9}\,\mathsf{P}(Z \ge x) \quad \forall x \in \mathbb{R}; \tag{1.3}$$

a multivariate analogue of (1.3) was also obtained there.

Later it was realized (Pinelis [22]) that it is possible to extract (1.3) from (1.1) for all $f \in \mathcal{F}^{(3)}$ because the tail function of the normal distribution is log-concave. The following is a special case of Theorem 4 of Pinelis [23]; see also Theorem 3.11 of Pinelis [22].

**Theorem 1.3.** *Suppose that $\alpha > 0$, $\xi$ and $\eta$ are real-valued r.v.'s, and the tail function $u \mapsto \mathsf{P}(\eta \ge u)$ is log-concave on $\mathbb{R}$. Then the comparison inequality*

$$\mathsf{E} f(\xi) \le \mathsf{E} f(\eta) \quad \forall f \in \mathcal{F}_+^{(\alpha)} \tag{1.4}$$

*implies*

$$\mathsf{P}(\xi \ge x) \le c_\alpha\,\mathsf{P}(\eta \ge x) \quad \forall x \in \mathbb{R}, \tag{1.5}$$

*where*

$$c_\alpha := \Gamma(\alpha+1)(e/\alpha)^\alpha.$$

*Moreover, the constant factor $c_\alpha$ is the best possible one in (1.5).*

A similar result for the special case when $\alpha = 1$ is due to Kemperman and is contained in the book by Shorack and Wellner [26, pages 797–799].

**Remark 1.4.** As follows from [22, Remark 3.13], a useful point is that the requirement of the log-concavity of the tail function $q(u) := \mathsf{P}(\eta \ge u)$ in Theorem 1.3 can be relaxed by replacing $q$ with any (e.g., the least) log-concave majorant of $q$. However, then the optimality of the constant factor $c_\alpha$ is not guaranteed.

Note that $c_3 = 2e^3/9$, which is the constant factor in (1.3). Bobkov, Götze, and Houdre [4] discovered a simpler proof of a variant of (1.3) with a constant factor $12.0099\dots$ in place of $2e^3/9 = 4.4634\dots$. The value of the constant factor is obviously important in statistical applications. The upper bound in (1.3) improves Chernoff-Hoeffding's bound $e^{-x^2/2}$ in (1.2) for all $x > 1.3124\dots$. On the other hand, the bound in [4] does so only for all $x > 4.5903\dots$, when $\mathsf{P}(Z \ge x) < 2.22 \times 10^{-6}$. The proof in [4] was direct (rather than based on a moment comparison inequality of the form (1.4)). Of course, this does not imply that the direct methods are inferior. In fact, this author has certain ideas to combine the direct and indirect methods to further improve the constant factors.

A stronger, "discrete" version of (1.3) was obtained in [23, Theorem 5], as follows. Let $\eta_1, \dots, \eta_n$ be independent zero-mean r.v.'s such that $|\eta_i| \le 1$ almost surely (a.s.) for all $i$, and let $b_1, \dots, b_n$ be any real numbers such that $b_1^2 + \cdots + b_n^2 = n$. Then

$$\mathsf{P}(b_1\eta_1 + \cdots + b_n\eta_n \ge x) \le c_3\, \mathsf{P}(\varepsilon_1 + \cdots + \varepsilon_n \ge x) \tag{1.6}$$

for all values $x$ that are taken on by the r.v. $\varepsilon_1 + \cdots + \varepsilon_n$ with nonzero probability. Clearly, (1.3) follows from (1.6) by the central limit theorem.

In this paper, we provide new upper bounds on generalized moments and tails of real-valued (super)martingales. It is well known that such bounds can be used, in particular, to obtain concentration-type results; see e.g. [17, 18, 22].

## 2. Upper bounds on generalized moments and tails of (super)martingales

**Theorem 2.1.** *Let $(S_0, S_1, \dots)$ be a supermartingale relative to a nondecreasing sequence of $\sigma$-algebras $H_{\le 0}, H_{\le 1}, \dots$, with $S_0 \le 0$ a.s. and differences $X_i := S_i - S_{i-1}$, $i = 1, 2, \dots$. Suppose that for every $i = 1, 2, \dots$ there exist non-random constants $d_i > 0$ and $\sigma_i > 0$ such that*

$$X_i \le d_i \quad \textit{and} \tag{2.1}$$

$$\mathsf{Var}(X_i | H_{\le i-1}) \le \sigma_i^2 \tag{2.2}$$

*a.s. Then, for all $n = 1, 2, \dots$,*

$$\mathsf{E} f(S_n) \le \mathsf{E} f(T_n) \quad \forall f \in \mathcal{F}_+^{(2)}, \quad \textit{where} \tag{2.3}$$

$$T_n := Z_1 + \cdots + Z_n \tag{2.4}$$

*and $Z_1, \dots, Z_n$ are independent r.v.'s such that each $Z_i$ takes on only two values, one of which is $d_i$, and satisfies the conditions*

$$\mathsf{E} Z_i = 0 \quad \textit{and} \quad \mathsf{Var} Z_i = \sigma_i^2; \quad \textit{that is,}$$

$$\mathsf{P}(Z_i = d_i) = \frac{\sigma_i^2}{d_i^2 + \sigma_i^2} \quad \textit{and} \quad \mathsf{P}\left(Z_i = -\frac{\sigma_i^2}{d_i}\right) = \frac{d_i^2}{d_i^2 + \sigma_i^2}.$$

Let us explain in general terms how such a result as Theorem 2.1 can be proved. First, it is not difficult to reduce Theorem 2.1 to the case when $(S_0, S_1, \dots)$ is a martingale with $S_0 = 0$ a.s. Then it is not difficult to reduce the situation to the case of one random summand $X$, so that the problem becomes to find – for any given $t \in \mathbb{R}$, $d > 0$, and $\sigma > 0$ – the maximum of $\mathsf{E} f_t(X)$ subject to the restrictions $X \le d$ a.s., $\mathsf{E}X = 0$, and $\mathsf{E}X^2 \le \sigma^2$, where $f_t(x) := (x-t)_+^2$. Then, using arguments going back to Chebyshev and Hoeffding [11, 13, 14, 15], one sees that here an extremal r.v. $X$ takes on at most three distinct values with nonzero probability. In fact, because of a special relation between the objective-function $f_t(x) = (x-t)_+^2$ and the restriction-functions 1, $x$, and $x^2$, one can see that an extremal r.v. $X$ takes on only two distinct values with nonzero probability; this makes the distribution of an extremal r.v. $X$ uniquely determined by $d$ and $\sigma$; in particular, the extremal distribution does not depend on the value of the parameter $t$ of the objective-function $f_t$. Thus the result follows.

An alternative approach is based on duality [22, Eq. (4)], and such an approach is actually used in the proof of Theorem 2.1 given in Section 3 below. According to the latter approach, one should search for a tight upper bound on the objective-function $f_t(x)$ over all $x \le d$; this upper bound must be of the form $Ax^2 + Bx + C$ (which is a linear combination of the restriction-functions 1, $x$, and $x^2$), so that $f_t(x) \le Ax^2 + Bx + C$ for all $x \le d$. It is not difficult to see that, for such a tight upper bound, the equality $f_t(x) = Ax^2 + Bx + C$ is attained for at most two distinct values of $x \in (-\infty, d]$. This again implies that an extremal r.v. $X$ takes on only two distinct values with nonzero probability, whence the result. (See the mentioned proof for details.)

**Remark 2.2.** In the case when $(S_i)$ is a martingale, Theorem 2.1 is a result of Bentkus [1, 3], who used essentially the "duality" approach. Moreover, using Schur convexity arguments similar to those in Eaton [5], he also showed that, in the case $d_i \equiv d$, for every $f \in \mathcal{F}_+^{(2)}$ the right-hand side $\mathsf{E} f(T_n)$ of inequality (2.3) is maximized – under the condition $\sigma_1^2 + \dots + \sigma_n^2 = n\sigma^2 = \text{const}$ – when the $Z_i$'s are identically distributed; this fact was earlier established by Hoeffding [12, (2.10)] for $f(x) \equiv e^{\lambda x}$ with $\lambda > 0$. Finally, for such i.i.d. $Z_i$'s with $d_i \equiv d$ and $\sigma_i \equiv \sigma$, Bentkus used the method given in [22] (cf. Theorem 1.3 and Remark 1.4 above) to extract the upper bound of the form

$$\mathsf{P}(S_n \ge y) \le c_2 \mathsf{P}^{\mathsf{LC}}(T_n \ge y) \quad \forall y \in \mathbb{R}, \tag{2.5}$$

where the function $y \mapsto \mathsf{P}^{\mathsf{LC}}(T_n \ge y)$ is the least log-concave majorant of the tail function $y \mapsto \mathsf{P}(T_n \ge y)$ on $\mathbb{R}$. Note that $c_2 = e^2/2 = 3.694\dots$ .

Note also that the distribution of r.v. $T_n$ in (2.5) is a shifted and re-scaled binomial distribution, concentrated on the lattice, say $L_{n,d,h}$, of all points of the form $nd + kh$ $(k \in \mathbb{Z})$, where

$$h := d + \sigma^2/d.$$

Here and henceforth, we assume that $d_i \equiv d$, unless indicated otherwise.

Since the tail function of the binomial distribution is log-concave on $\mathbb{Z}$ (see e.g. [23, Remark 13]), one has $\mathsf{P}^{\mathsf{LC}}(T_n \ge y) = \mathsf{P}(T_n \ge y)$ for all $y$ in the lattice $L_{n,d,h}$.

Inequality (2.5) can be significantly improved. Let the function $y \mapsto \mathsf{P}^{\mathsf{Lin}}(T_n \ge y)$ denote the linear interpolation of the function $y \mapsto \mathsf{P}(T_n \ge y)$ over the lattice $L_{n,d,h}$, so that

$$\mathsf{P}^{\mathsf{Lin}}(T_n \ge y) = (1-\gamma)\mathsf{P}(T_n \ge nd + kh) + \gamma \mathsf{P}(T_n \ge nd + (k+1)h)$$

if $\gamma := (y - nd - kh)/h \in [0,1]$ for some $k \in \mathbb{Z}$. Let then the function $y \mapsto \mathsf{P}^{\mathsf{Lin,LC}}(T_n \geq y)$ denote the least log-concave majorant of the function $y \mapsto \mathsf{P}^{\mathsf{Lin}}(T_n \geq y)$ on $\mathbb{R}$.

**Theorem 2.3.** *Suppose that the conditions of Theorem 2.1 hold with $d_i \equiv d$. As before, let $\sigma^2 := \frac{1}{n}(\sigma_1^2 + \cdots + \sigma_n^2)$ and $h := d + \sigma^2/d$. Then*

$$\mathsf{P}(S_n \geq y) \leq c_2 \mathsf{P}^{\mathsf{Lin,LC}}(T_n \geq y + \tfrac{h}{2}) \quad \forall y \in \mathbb{R}. \tag{2.6}$$

Because the tail function $y \mapsto \mathsf{P}(T_n \geq y)$ decreases very rapidly, the shift $\frac{h}{2}$ in $\mathsf{P}^{\mathsf{Lin,LC}}(T_n \geq y + \frac{h}{2})$ generally provides quite a substantial improvement. As will be shown later in this paper (Proposition 2.7), in almost all practically important cases the bound (2.6) is better, or even much better, than (2.5).

In Subsection 2.1 (Proposition 2.10), we will also provide an explicit expression for $\mathsf{P}^{\mathsf{Lin,LC}}(T_n \geq y + \frac{h}{2})$.

That $(S_0, S_1, \dots)$ in Theorems 2.1 and 2.3 is allowed to be a supermartingale (rather than only a martingale) makes it convenient to use the simple but powerful truncation tool. (Such a tool was used, for example, in [20] to prove limit theorems for large deviation probabilities in Banach spaces based only on precise enough probability inequalities and without using Cramér's transform, the standard device in the theory of large deviations.) Thus, for instance, one immediately has the following corollary of Theorem 2.3.

**Corollary 2.4.** *Suppose that all conditions of Theorem 2.1 hold except possibly for condition* (2.1)*. Then for all $y \in \mathbb{R}$ and $d > 0$*

$$\mathsf{P}(S_n \geq y) \leq \mathsf{P}\big(\max_{1 \leq i \leq n} X_i \geq d\big) + c_2\, \mathsf{P}^{\mathsf{Lin,LC}}(T_n \geq y + \tfrac{h}{2}) \tag{2.7}$$

$$\leq \sum_{1 \leq i \leq n} \mathsf{P}\,(X_i \geq d) + c_2\, \mathsf{P}^{\mathsf{Lin,LC}}(T_n \geq y + \tfrac{h}{2}). \tag{2.8}$$

These bounds are much more precise than the exponential bounds in [8, 9, 19].

**Remark 2.5.** By the Doob inequality, inequality (2.6) holds for the maximum, $M_n := \max_{0 \leq k \leq n} S_k$, in place of $S_n$. This follows because (i) all functions of class $\mathcal{F}_+^{(2)}$ are convex and (ii) in view of Lemma 3.1 on page 41, one may assume without loss of generality that $(S_i)$ is a martingale. Similarly, inequalities (2.7) and (2.8) hold for $M_n$ in place of $S_n$.

In a similar manner, under conditions (2.1)–(2.2) and with

$$b := \sqrt{b_1^2 + \cdots + b_n^2}, \quad \text{where} \quad b_i := \max(d_i, \sigma_i), \tag{2.9}$$

Bentkus [2, 3] obtained the following extensions of inequalities (1.6) and (1.3), respectively:

$$\mathsf{P}(S_n \geq y) \leq c_3 \mathsf{P}^{\mathsf{LC}}\Big(\sum_{i=1}^n \varepsilon_i \geq \frac{y\sqrt{n}}{b}\Big) \quad \forall y \in \mathbb{R} \quad \text{and} \tag{2.10}$$

$$\mathsf{P}(S_n \geq y) \leq c_3 \mathsf{P}(Z \geq y/b) \quad \forall y \in \mathbb{R}.$$

The upper bound in (2.10) can be improved in a similar manner, as follows.

**Proposition 2.6.** *Under conditions* (2.1)*,* (2.2)*, and* (2.9)*,*

$$\mathsf{P}(S_n \geq y) \leq c_3 \mathsf{P}^{\mathsf{Lin,LC}}\Big(\sum_{i=1}^n \varepsilon_i \geq 1 + \frac{y\sqrt{n}}{b}\Big) \quad \forall y \in \mathbb{R}. \tag{2.11}$$

In (2.10) and (2.11), the $d_i$'s are allowed to differ from one another.

Note that the expression $\mathsf{P}^{\mathsf{Lin,LC}}\big(\sum_{i=1}^n \varepsilon_i \ge 1 + \frac{y\sqrt{n}}{b}\big) = \mathsf{P}^{\mathsf{Lin,LC}}\big(\sum_{i=1}^n \frac{b}{\sqrt{n}}\varepsilon_i \ge y+\frac{b}{\sqrt{n}}\big)$ is a special case of the expression $\mathsf{P}^{\mathsf{Lin,LC}}(T_n \ge y+\frac{h}{2})$, with $\frac{h}{2} = d = \sigma = \frac{b}{\sqrt{n}}$.

From the "right-tail" bounds stated above, "left-tail" and "two-tail" ones immediately follow. For instance, if condition $X_i \le d_i$ in Theorem 2.1 is replaced by $|X_i| \le d_i$, then inequality (2.3) holds with $\mathcal{F}^{(2)}$ in place of $\mathcal{F}_+^{(2)}$ provided that $(S_0, S_1, \dots)$ is a martingale with $S_0 = 0$ a.s.

In order to present an explicit formula for the upper bound in (2.6) and compare it with the upper bound in (2.5), it is convenient to rescale the r.v. $T_n$, taking on values in the lattice $L_{n,d,h}$ of all points of the form $nd + kh$ ($k \in \mathbb{Z}$), so that the rescaled r.v., say

$$B_n := \frac{q}{d}T_n + np, \quad \text{where} \quad p := \frac{\sigma^2}{d^2+\sigma^2} \quad \text{and} \quad q := 1-p = \frac{d^2}{d^2+\sigma^2}, \tag{2.12}$$

is binomial with parameters $n$ and $p$. Then for all $y \in \mathbb{R}$, with

$$x := \frac{q}{d}y + np, \tag{2.13}$$

one has

$$\begin{aligned}
\mathsf{P}(T_n \ge y) &= \mathsf{P}(B_n \ge x) =: Q_n(x),\\
\mathsf{P}^{\mathsf{LC}}(T_n \ge y) &= \mathsf{P}^{\mathsf{LC}}(B_n \ge x) =: Q_n^{\mathsf{LC}}(x),\\
\mathsf{P}^{\mathsf{Lin}}(T_n \ge y) &= \mathsf{P}^{\mathsf{Lin}}(B_n \ge x) =: Q_n^{\mathsf{Lin}}(x),\\
\mathsf{P}^{\mathsf{Lin,LC}}(T_n \ge y) &= \mathsf{P}^{\mathsf{Lin,LC}}(B_n \ge x) =: Q_n^{\mathsf{Lin,LC}}(x), \quad \text{so that}\\
\mathsf{P}^{\mathsf{Lin,LC}}(T_n \ge y+\tfrac{h}{2}) &= Q_n^{\mathsf{Lin,LC}}(x+\tfrac{1}{2}). \qquad (2.14)
\end{aligned}$$

Here the function $x \mapsto \mathsf{P}^{\mathsf{Lin,LC}}(B_n \ge x)$ is defined as the least log-concave majorant on $\mathbb{R}$ of the function $x \mapsto \mathsf{P}^{\mathsf{Lin}}(B_n \ge x)$, which is in turn defined as the linear interpolation of the tail function $x \mapsto \mathsf{P}(B_n \ge x)$ over the lattice $\mathbb{Z}$. Similarly, the function $x \mapsto \mathsf{P}^{\mathsf{LC}}(B_n \ge x)$ is defined as the least log-concave majorant of the function $x \mapsto \mathsf{P}(B_n \ge x)$ on $\mathbb{R}$.

Note also that $B_n \in [0, n]$ a.s.

Now one is ready to state the following comparison between the upper bounds in (2.5) and (2.6).

**Proposition 2.7.** *Here relation* (2.13) *between $y$ and $x$ is assumed.*

**(i)** *Equation*

$$\ln\frac{1-u}{-\ln u} - 1 - \frac{1}{2}\,\frac{(1+u)\ln u}{1-u} = 0 \tag{2.15}$$

*in $u \in (0,1)$ has exactly one solution, $u = u_* := 0.00505778\ldots$ .*

**(ii)** *The upper bound in* (2.6) *is no greater than that in* (2.5) *for all $x \le j_{**}$, where*

$$j_{**} := \left\lfloor \frac{n - u_{**}\frac{q}{p}}{1 + u_{**}\frac{q}{p}} \right\rfloor \quad \textit{and} \tag{2.16}$$

$$u_{**} := \frac{u_*}{1-u_*} = 0.00508349\ldots; \tag{2.17}$$

*since $u_{**}$ is small, $j_{**}$ is rather close to $n$ unless $\frac{q}{p}$ is large.*

**(iii)** *Moreover, the upper bound in* (2.6) *is no greater than that in* (2.5) *for all* $x \le n$ *provided that*

$$n \le \frac{p}{q}\,\frac{1}{u_{**}} = \frac{p}{q}\,196.714\ldots\ . \tag{2.18}$$

*In particular,* (2.6) *works better than* (2.5) *for all* $x \le n$ *if* $n \le 196$ *and* $p \ge \frac{1}{2}$.

Since $\sum_{i=1}^{n} \varepsilon_i$ is a special case of $T_n$ (with $d = \sigma = 1$), one immediately obtains the following corollary to Proposition 2.7.

**Corollary 2.8.** *The upper bound in* (2.11) *is no greater than that in* (2.10) *for all* $x \le j_{**}$*, where* $x := \frac{y\sqrt{n}}{2b} + \frac{n}{2}$ *and* $j_{**} := \left\lfloor \frac{n-u_{**}}{1+u_{**}} \right\rfloor$. *Moreover, the upper bound in* (2.11) *is no greater than that in* (2.10) *for all* $x \le n$ *provided that* $n \le 196$.

Of course, the restriction $x \le j_{**}$ (even though very weak) is only sufficient, but not necessary for the upper bound in (2.6) to be no greater than that in (2.5).

Moreover, (2.6) may work very well even when $p$ is small. For example, in Figure 1 one can see the graph of the ratio

$$r(x) := \frac{Q_n^{\mathsf{Lin,LC}}\left(x + \frac{1}{2}\right)}{Q_n^{\mathsf{LC}}\,(x)} = \frac{c_2\,\mathsf{P}^{\mathsf{Lin,LC}}(T_n \ge y + \frac{h}{2})}{c_2\,\mathsf{P}^{\mathsf{LC}}(T_n \ge y)}$$

of the "new" upper bound – that in (2.6), to the "old" one – that in (2.5), for $n = 30$ and $p = \frac{3}{100}$. In the same figure, one can also see the graph of the new upper bound $q(x) := \min\left(1, c_2\,Q_n^{\mathsf{Lin,LC}}\left(x + \frac{1}{2}\right)\right)$, which is a very rapidly decreasing tail function. Proposition 2.7(ii) guaranteed that, for these $n$ and $p$, the new upper bound will be an improvement of the old one (that is, one will have $r(x) \le 1$) at least for all $x \in (-\infty, j_{**}] = (-\infty, 25]$, which is rather close to what one can see in the picture. Note that, by Theorem 2.3, for $x \ge 25$ one has $\mathsf{P}(S_n \ge y) \le c_2 Q_n^{\mathsf{Lin,LC}}\left(25 + \frac{1}{2}\right) \approx 3.44 \times 10^{-33}$. Thus, the new upper bound is not an improvement only if the tail probability $\mathsf{P}(S_n \ge y)$ is less than $3.45 \times 10^{-33}$. On the other hand, for instance, $r(4) \approx 0.58$; that is, the new upper bound is approximately a 42% improvement of the old one for $x = 4$; at that, the new upper bound is $\approx 0.026 = 2.6\%$, a value quite in a common range in statistical practice.

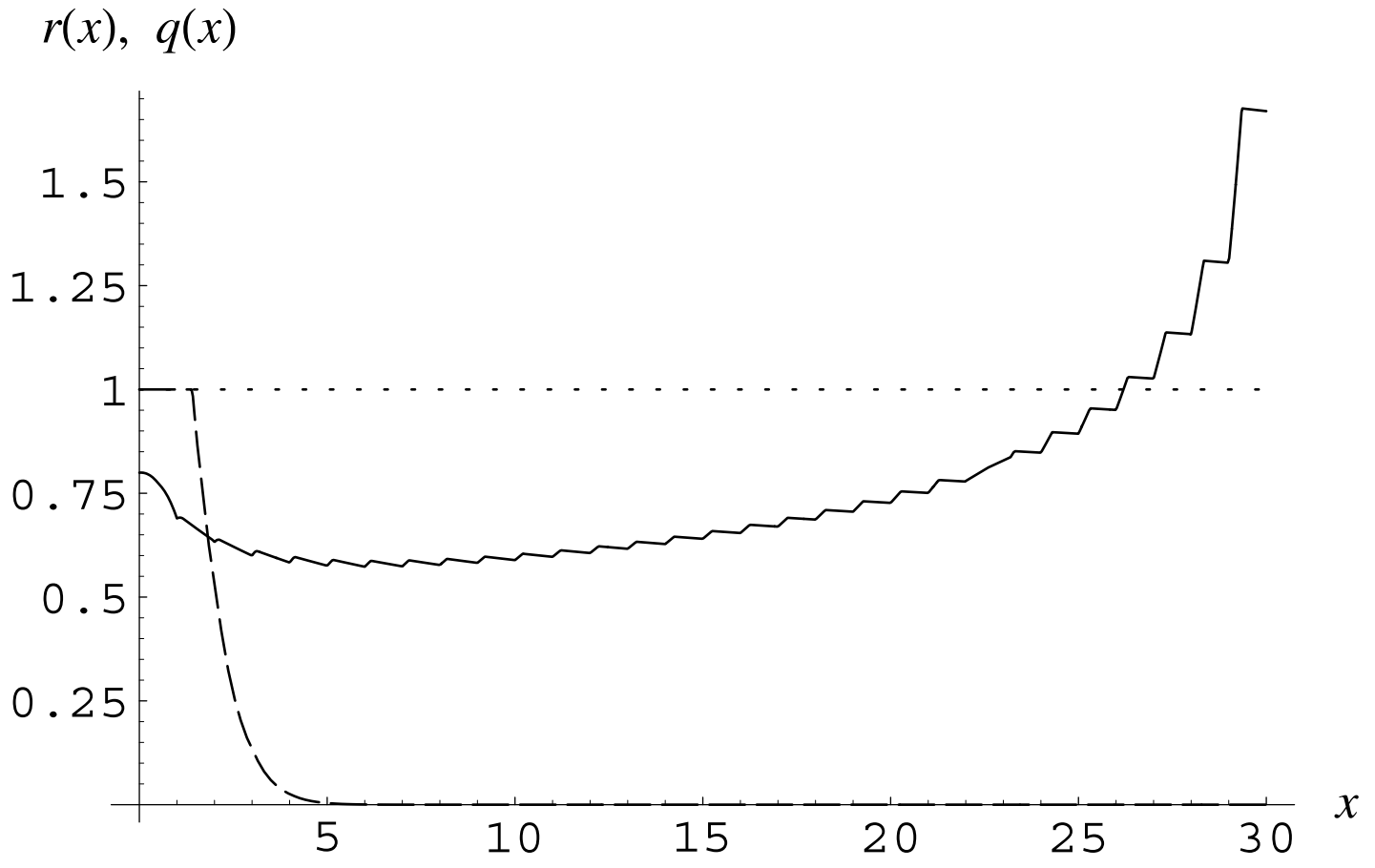

Fig 1. $r(x)$, *solid;* $q(x) = \min\left(1, c_2\,Q_n^{\mathsf{Lin,LC}}\left(x + \frac{1}{2}\right)\right)$, *long dashes.*

The upper bound in (2.6), $c_2\,\mathsf{P}^{\mathsf{Lin,LC}}(T_n \ge y + \frac{h}{2}) = c_2\,Q_n^{\mathsf{Lin,LC}}\left(x + \frac12\right)$, is rather simple to compute, as described in Proposition 2.10, given in a separate subsection, Subsection 2.1. An underlying reason for this simplicity is that the "discrete" tail function, $\mathbb{Z} \ni j \mapsto Q_n(j)$, of the binomial distribution is log-concave [23, Remark 13]; therefore, it turns out (by Propositions 2.10 and 2.9) that the value $Q_n^{\mathsf{Lin,LC}}\left(x + \frac12\right)$ can be computed locally: for a certain function $\mathcal{Q}$ which depends only on its 7 arguments, one has $Q_n^{\mathsf{Lin,LC}}\left(x + \frac12\right) = \mathcal{Q}(x, n, j_*, q_{k-1}, q_k, q_{k+1}, q_{k+2})$ for all $x \in \mathbb{R}$ and $p \in (0,1)$, where $j_* := \lfloor (n+1)p \rfloor + 1$, $k := \lfloor x \rfloor$, and $q_j := Q_n(j)$ for all $j \in \mathbb{Z}$.

## 2.1. How to compute the expression $Q_n^{\mathsf{Lin,LC}}\left(x + \frac12\right)$ in (2.14)

First here, we need to introduce some notation.

For each $j = 1, \dots, n$ such that $j > (n+1)p$, let

$$x_j := j - \frac12 + \frac{q_j}{p_j} + \frac{\frac{q_j}{p_{j-1}} - \frac{q_j}{p_j}}{\ln \frac{p_{j-1}}{p_j}}, \tag{2.19}$$

$$y_j := j - \frac12 + \frac{q_j}{p_{j-1}} + \frac{\frac{q_j}{p_{j-1}} - \frac{q_j}{p_j}}{\ln \frac{p_{j-1}}{p_j}} = x_j - q_j\left(\frac{1}{p_j} - \frac{1}{p_{j-1}}\right), \tag{2.20}$$

where, for all $j \in \mathbb{Z}$,
(2.21)

$$q_j := Q_n(j) = \mathsf{P}(B_n \ge j) \quad \text{and} \quad p_j := q_j - q_{j+1} = \mathsf{P}(B_n = j) = \binom{n}{j} p^j q^{n-j};$$

note that, for each $j = 0, \dots, n$,

$$p_{j-1} > p_j \iff j > (n+1)p \iff j \ge j_* := \lfloor (n+1)p \rfloor + 1, \tag{2.22}$$

so that, for all $j \in \mathbb{Z} \cap [j_*, n]$, one has $\ln \frac{p_{j-1}}{p_j} > 0$, and so, $x_j$ and $y_j$ are well defined, and $x_j > y_j$; note that $j_* \ge 1$. Let also

$$x_j := j + \tfrac12 \quad \text{and} \quad y_j := j - \tfrac12 \quad \text{for integer } j \ge n+1. \tag{2.23}$$

**Proposition 2.9.** *For all integer $j \ge j_*$, one has*

$$j - \tfrac32 < j - 1 < y_j \le j - \tfrac12 < x_j \le y_{j+1} \le j + \tfrac12;$$

*moreover, if $j \le n$, then $y_j < j - \frac12$ and $x_j < y_{j+1}$.*

By Proposition 2.9, for all integer $j \ge j_*$, the intervals

$$\delta_j := (y_j, x_j)$$

are non-empty, with the endpoints

$$y_j \in (j-1, j-\tfrac12] \subset (j-\tfrac32, j-\tfrac12] \quad \text{and} \quad x_j \in (j-\tfrac12, j+\tfrac12];$$

moreover, the intervals $\delta_j$ are strictly increasing in $j \ge j_*$: $\delta_j < \delta_{j+1}$ for all $j \ge j_*$, where we use the following convention for any two subsets $A$ and $B$ of $\mathbb{R}$:

$$A < B \overset{\text{def}}{\iff} (x < y \quad \forall x \in A\ \forall y \in B).$$

For all integer $j \ge j_*$ and all $x \in \delta_j$, introduce the interpolation expression

$$Q_n^{\mathsf{Interp}}(x; j) := Q_n^{\mathsf{Lin}}\left(y_j + \tfrac12\right)^{1-\delta} Q_n^{\mathsf{Lin}}\left(x_j + \tfrac12\right)^{\delta}, \quad \text{where } \delta := \frac{x - y_j}{x_j - y_j}. \tag{2.24}$$

**Proposition 2.10.** *For all real $x$,*

$$Q_n^{\mathsf{Lin,LC}}\left(x+\tfrac{1}{2}\right)=Q(x):=\begin{cases}Q_n^{\mathsf{Interp}}(x;j) & \text{if } \exists j\in\mathbb{Z}\cap[j_*,n]\ \ x\in\delta_j,\\ Q_n^{\mathsf{Lin}}\left(x+\tfrac{1}{2}\right) & \text{otherwise.}\end{cases}\tag{2.25}$$

A few comments are in order here:

- the function $Q$ in (2.25) is well defined, because the $\delta_j$'s are pairwise disjoint;
- $Q_n^{\mathsf{Lin,LC}}\left(x+\frac{1}{2}\right)$ is easy to compute by (2.25) because, in view of Proposition 2.9, the condition $x\in\delta_j$ for $j\in\mathbb{Z}\cap[j_*,n]$ implies that $j$ equals either $\lfloor x\rfloor$ or $\lfloor x\rfloor+1$. In particular, $Q_n^{\mathsf{Lin,LC}}\left(x+\frac{1}{2}\right)=Q_n^{\mathsf{Lin}}\left(x+\frac{1}{2}\right)=0$ for all $x\geq n+\frac{1}{2}$.

## 3. Proofs

*Proof of Proposition 1.1.* The "only if" part follows because $f\in\mathcal{F}^{(\alpha)}$ for a natural $\alpha$ implies that the (right) derivative $f'(u)=\alpha\int_{-\infty}^{\infty}(u-t)_+^{\alpha-1}\,\mu(\mathrm{d}t)$. Vice versa, if $f$ satisfies the conditions listed after "iff" in the statement of Proposition 1.1, then one can use the Fubini theorem repeatedly to see that for all real $u$

$$f(u)=\int_{-\infty}^{u}f'(t)\,\mathrm{d}t=\int_{-\infty}^{u}\mathrm{d}t\int_{-\infty}^{t}f''(s)\,\mathrm{d}s=\int_{-\infty}^{u}(u-s)f''(s)\,\mathrm{d}s=\dots$$
$$=\int_{-\infty}^{u}\frac{(u-s)^{\alpha-1}}{(\alpha-1)!}f^{(\alpha)}(s)\,\mathrm{d}s=\int_{-\infty}^{u}\frac{(u-s)^{\alpha}}{\alpha!}\,\mathrm{d}f^{(\alpha)}(s)=\int_{-\infty}^{\infty}(u-s)_+^{\alpha}\,\mu(\mathrm{d}s),$$

where $f^{(\alpha)}$ is the (nondecreasing) right derivative of the convex function $f^{(\alpha-1)}$ and $\mu(\mathrm{d}s):=\mathrm{d}f^{(\alpha)}(s)/\alpha!$. □

Theorem 2.1 can be rather easily reduced to the case when $(S_n)$ is a martingale. This is implied by the following lemma.

**Lemma 3.1.** *Let $(S_i)$ be a supermartingale as in Theorem 2.1, so that conditions (2.1) and (2.2) are satisfied. For $i=1,2,\dots,$ let*

$$\tilde{X}_i:=(1-\gamma_{i-1})X_i+\gamma_{i-1}d_i,\quad \text{where}\quad \gamma_{i-1}:=\frac{\mathsf{E}_{i-1}X_i}{\mathsf{E}_{i-1}X_i-d_i};\tag{3.1}$$

*$\mathsf{E}_j$ and $\mathsf{Var}_j$ denote, respectively, the conditional expectation and variance given $X_1,\dots,X_j$. Then $\tilde{X}_i$ is $H_{\leq i}$-measurable,*

$$X_i\leq\tilde{X}_i\leq d_i,\quad \mathsf{E}_{i-1}\tilde{X}_i=0,\quad \text{and}\quad \mathsf{Var}_{i-1}\tilde{X}_i\leq\mathsf{Var}_{i-1}X_i\leq\sigma_i^2\quad \text{a.s.}$$

*Proof.* The conditions $d_i>0$ and $\mathsf{E}_{i-1}X_i\leq 0$ imply that that $\gamma_{i-1}\in[0,1)$. Now (3.1) and the inequality $X_i\leq d_i$ yield $X_i\leq\tilde{X}_i\leq d_i$. Moreover, (3.1) yields $\mathsf{E}_{i-1}\tilde{X}_i=(1-\gamma_{i-1})\mathsf{E}_{i-1}X_i+\gamma_{i-1}d_i=0$ and $\mathsf{Var}_{i-1}\tilde{X}_i=(1-\gamma_i)^2\,\mathsf{Var}_{i-1}X_i\leq\mathsf{Var}_{i-1}X_i$. □

Theorem 2.1 is mainly based on the following lemma, which also appeared as [1, (12)] and [3, Lemma 4.4(i), with condition $\mathsf{E}X=0$ missing there].

**Lemma 3.2.** *Let $X$ be a r.v. such that $\mathsf{E}X=0$, $\mathsf{Var}X\leq\sigma^2$, and $X\leq d$ a.s. for some non-random constants $\sigma>0$ and $d>0$. Let $a:=\sigma^2/d^2$, and let $X_a$ be a r.v. taking on values $-a$ and $1$ with probabilities $\frac{1}{1+a}$ and $\frac{a}{1+a}$, respectively. Then*

$$\mathsf{E}f(X)\leq\mathsf{E}f(d\cdot X_a)\quad\forall f\in\mathcal{F}_+^{(2)}.\tag{3.2}$$

*Proof.* (Given here for the reader's convenience and because it is short and simple.) By homogeneity, one may assume that $d = 1$ (otherwise, rewrite the lemma in terms of $X/d$ in place of $X$). Note that $X_a \le 1$ with probability 1, $\mathsf{E}X_a = 0$, and $\mathsf{E}X_a^2 = \sigma^2$. Let here $f_t(x) := (x-t)_+^2$ and

$$h_t(x) := \frac{(1-t)_+^2}{(1-t_a)^2}\,(x-t_a)^2, \quad \text{where} \quad t_a := \min(t, -a).$$

Then it is easy to check (by considering the cases $t \ge 1$, $-a \le t \le 1$, and $t \le -a$) that $f_t(x) \le h_t(x)$ for all $x \le 1$, and $f_t(x) = h_t(x)$ for $x \in \{1, -a\}$. Therefore, $\mathsf{E}f_t(X) \le \mathsf{E}h_t(X) \le \mathsf{E}h_t(X_a) = \mathsf{E}f_t(X_a)$ (the second inequality here follows because $h_t(x)$ is a quadratic polynomial in $x$ with a nonnegative coefficient of $x^2$, while $\mathsf{E}X = \mathsf{E}X_a$ and $\mathsf{E}X^2 \le \mathsf{E}X_a^2$). Now the lemma follows by the definition of the class $\mathcal{F}_+^{(\alpha)}$ and the Fubini theorem. □

*Proof of Theorem 2.1.* This proof is based in a standard manner on Lemma 3.2, using also Lemma 3.1. Indeed, by Lemma 3.1, one may assume that $\mathsf{E}_{i-1}X_i = 0$ for all $i$. Let $Z_1, \dots, Z_n$ be r.v.'s as in the statement of Theorem 2.1, which are also independent of the $X_i$'s, and let

$$R_i := X_1 + \cdots + X_i + Z_{i+1} + \cdots + Z_n.$$

Let $\tilde{\mathsf{E}}_i$ denote the conditional expectation given $X_1, \dots, X_{i-1}, Z_{i+1}, \dots, Z_n$. Note that, for all $i = 1, \dots, n$, $\tilde{\mathsf{E}}_i X_i = \mathsf{E}_{i-1}X_i = 0$ and $\tilde{\mathsf{E}}_i X_i^2 = \mathsf{E}_{i-1}X_i^2$; moreover, $R_i - X_i = X_1 + \cdots + X_{i-1} + Z_{i+1} + \cdots + Z_n$ is a function of $X_1, \dots, X_{i-1}$, $Z_{i+1}, \dots, Z_n$. Hence, by Lemma 3.2, for any $f \in \mathcal{F}_+^{(2)}$, $\tilde{f}_i(x) := f(R_i - X_i + x)$, and all $i = 1, \dots, n$, one has $\tilde{\mathsf{E}}_i f(R_i) = \tilde{\mathsf{E}}_i \tilde{f}_i(X_i) \le \tilde{\mathsf{E}}_i \tilde{f}_i(Z_i) = \tilde{\mathsf{E}}_i f(R_{i-1})$, whence $\mathsf{E}f(S_n) \le \mathsf{E}f(R_n) \le \mathsf{E}f(R_0) = \mathsf{E}f(T_n)$ (the first inequality here follows because $S_0 \le 0$ a.s. and any function $f$ in $\mathcal{F}_+^{(2)}$ is nondecreasing). □

*Proof of Theorem 2.3.* In view of Theorem 2.1 and Remark 2.2, one has

$$\mathsf{E}f(\tilde{S}_n) \le \mathsf{E}f(B_n) \quad \forall f \in \mathcal{F}_+^{(2)}, \tag{3.3}$$

where $\tilde{S}_n := \frac{q}{d}S_n + np$ and $B_n$ is defined by (2.12). Let $U$ be a r.v., which is independent of $B_n$ and uniformly distributed between $-\frac{1}{2}$ and $\frac{1}{2}$. Then, by Jensen's inequality, $\mathsf{E}f(B_n) \le \mathsf{E}f(B_n + U)$ for all convex functions $f$, whence

$$\mathsf{E}f(\tilde{S}_n) \le \mathsf{E}f(B_n + U) \quad \forall f \in \mathcal{F}_+^{(2)}.$$

Observe that the density function of $B_n + U$ is $x \mapsto \sum_{j=0}^n p_j\, \mathbf{I}\{|x - j| < \frac{1}{2}\}$ (where the $p_j$'s are given by (2.21)), and so, the tail function of $B_n + U$ is given by the formula $\mathsf{P}(B_n + U \ge x) = Q_n^{\mathsf{Lin}}\big(x + \frac{1}{2}\big)$ $\forall x \in \mathbb{R}$. Now Theorem 2.3 follows by Theorem 1.3, Remark 1.4, and (2.14). □

*Proof of Proposition 2.6.* This proof is quite similar to the proof of Theorem 2.3. Instead of (3.3), here one uses inequality $\mathsf{E}f(\tilde{S}_n) \le \mathsf{E}f(B_n)$ $\forall f \in \mathcal{F}_+^{(3)}$, where $\tilde{S}_n := \frac{\sqrt{n}}{2b}S_n + \frac{n}{2}$. This latter inequality follows from (2.3) (with $d_i$ and $\sigma_i$ each replaced by $b_i = \max(d_i, \sigma_i)$) and the first one of the inequalities (1.1) $\forall f \in \mathcal{F}_+^{(3)}$ (recall Remark 1.2), taking also into account the inclusion $\mathcal{F}_+^{(3)} \subseteq \mathcal{F}_+^{(2)}$ (which follows e.g. from [23, Proposition 1(ii)]). □

In the following two propositions, which are immediate corollaries of results of [25] and [24], it is assumed that $f$ and $g$ are differentiable functions on an interval $(a,b) \subseteq (-\infty,\infty)$, and each of the functions $g$ and $g'$ is nonzero and does not change sign on $(a,b)$; also, $r := f/g$ and $\rho := f'/g'$.

**Proposition 3.1.** (Cf. [25, Proposition 1].) *Suppose that* $f(a+) = g(a+) = 0$ *or* $f(b-) = g(b-) = 0$. *Suppose also that* $\rho \nearrow$ *or* $\searrow$*; that is,* $\rho$ *is increasing or decreasing on* $(a,b)$. *Then* $r \nearrow$ *or* $\searrow$*, respectively.*

**Proposition 3.2.** (Cf. [24, Proposition 1.9]; the condition that $f(a+) = g(a+) = 0$ or $f(b-) = g(b-) = 0$ is not assumed here.) *If* $\rho \nearrow$ *or* $\searrow$ *on* $(a,b)$, *then* $r \nearrow$ *or* $\searrow$ *or* $\nearrow\searrow$ *or* $\searrow\nearrow$ *on* $(a,b)$. *(Here, for instance, the pattern* $\nearrow\searrow$ *for* $\rho$ *on* $(a,b)$ *means that* $\rho \nearrow$ *on* $(a,c)$ *and* $\searrow$ *on* $(c,b)$ *for some* $c \in (a,b)$*; the pattern* $\searrow\nearrow$ *has a similar meaning.) It follows that, if* $\rho \nearrow$ *or* $\searrow$ *or* $\nearrow\searrow$ *or* $\searrow\nearrow$ *on* $(a,b)$, *then* $r \nearrow$ *or* $\searrow$ *or* $\nearrow\searrow$ *or* $\searrow\nearrow$ *or* $\nearrow\searrow\nearrow$ *or* $\searrow\nearrow\searrow$ *on* $(a,b)$.

**Lemma 3.3.** *Part* (i) *of Proposition 2.7 is true.*

*Proof of Lemma 3.3.* Let

$$h(u) := \ln\frac{1-u}{-\ln u} - 1 - \frac{1}{2}\,\frac{(1+u)\ln u}{1-u}, \tag{3.4}$$

the left-hand side of (2.15). Here and in rest of the proof of Lemma 3.3, it is assumed that $0 < u < 1$, unless specified otherwise. Then

$$h'(u) = r_1(u) := \frac{f_1(u)}{g_1(u)}, \tag{3.5}$$

where $f_1(u) := 2\ln^2 u + (\frac{1}{u} + 2 - 3u)\ln u + 2(u + \frac{1}{u}) - 4$ and $g_1(u) := -2(1-u)^2 \ln u$. Let next

$$r_2(u) := \frac{f_1'(u)}{g_1'(u)} = \frac{f_2(u)}{g_2(u)}, \tag{3.6}$$

where $f_2(u) := (\frac{3}{u} - \frac{1}{u^2})\ln u + \frac{1}{u} - \frac{1}{u^2}$ and $g_2(u) := 4\ln u - \frac{2}{u} + 2$, and then

$$r_3(u) := \frac{f_2'(u)}{g_2'(u)} = \frac{f_3(u)}{g_3(u)}, \tag{3.7}$$

where $f_3(u) := (2-3u)\ln u + 1 + 2u$ and $g_3(u) := 2u(1+2u)$.

One has $\frac{f_3''(u)}{g_3''(u)} = -\frac{1}{8}(\frac{2}{u^2} + \frac{3}{u})$, which is increasing; moreover, $\frac{\mathrm{d}}{\mathrm{d}u}\frac{f_3'(u)}{g_3'(u)}$ tends to $-\infty < 0$ and $-29/50 < 0$ as $u \downarrow 0$ and $u \uparrow 1$, respectively. Hence, by Proposition 3.2, $\frac{f_3'(u)}{g_3'(u)}$ is decreasing (in $u \in (0,1)$).

Next, by (3.7), $r_3'(0+) = \infty > 0$ and $r_3'(1-) = -2/3 < 0$. Hence, by Proposition 3.2, $r_3 \nearrow\searrow$ (on $(0,1)$).

By (3.6), $r_2'(0+) = \infty > 0$ and $f_2(1) = g_2(1) = 0$. Hence, by Propositions 3.1 and 3.2, $r_2 \nearrow\searrow$ (on $(0,1)$).

By (3.5), $r_1'(0+) = \infty > 0$ and $r_1'(1-) = -1/4 < 0$. Hence, by Proposition 3.2, $h' = r_1 \nearrow\searrow$ (on $(0,1)$). Moreover, $h'(0+) = -\infty < 0$ and $h'(1-) = \frac{1}{2} > 0$. Hence, for some $\beta \in (0,1)$, one has $h' < 0$ on $(0,\beta)$ and $h' > 0$ on $(\beta,1)$.

Hence, $h \searrow\nearrow$ on $(0,1)$. Moreover, $h(0+) = \infty$ and $h(1-) = 0$. It follows that the equation $h(u) = 0$ has a unique root $u = u_* \in (0,1)$.

Now Lemma 3.3 follows by (3.4). ☐

**Lemma 3.4.** *For all* $x \le j_{**} + 1$,

$$Q_n^{\mathsf{Lin}}\left(x + \tfrac{1}{2}\right) \le Q_n^{\mathsf{LC}}(x). \tag{3.8}$$

*Proof of Lemma 3.4.* For any given $x \le j_{**} + 1$, let

$$j := j_x := \lfloor x \rfloor \quad \text{and} \quad k := k_x := \lfloor x + \tfrac{1}{2} \rfloor, \quad \text{so that}$$

$$\begin{gathered} j \le x < j+1, \quad k - \tfrac{1}{2} \le x < k + \tfrac{1}{2}, \\ Q_n^{\mathsf{LC}}(x) = q_j^{1-\delta} q_{j+1}^{\delta}, \quad Q_n^{\mathsf{Lin}}\left(x + \tfrac{1}{2}\right) = (1-\gamma) q_k + \gamma\, q_{k+1}, \quad \text{where} \\ \delta := x - j \in [0,1) \quad \text{and} \quad \gamma := x + \tfrac{1}{2} - k \in [0,1). \end{gathered}$$

There are only three possible cases: $\delta = 0$, $\delta \in [\frac{1}{2}, 1)$, and $\delta \in (0, \frac{1}{2})$.

**Case 1:** $\delta = 0$. This case is simple. Indeed, here $k = j$ and $\gamma = \frac{1}{2}$, so that

$$Q_n^{\mathsf{LC}}(x) = q_j \ge \tfrac{1}{2}\, q_j + \tfrac{1}{2}\, q_{j+1} = Q_n^{\mathsf{Lin}}\left(x + \tfrac{1}{2}\right),$$

since $q_j$ is nonincreasing in $j$.

**Case 2:** $\delta \in [\frac{1}{2}, 1)$. This case is simple as well. Indeed, here $k = j+1$, so that

$$Q_n^{\mathsf{LC}}(x) \ge q_{j+1} \ge (1-\gamma)\, q_{j+1} + \gamma\, q_{j+2} = Q_n^{\mathsf{Lin}}\left(x + \tfrac{1}{2}\right).$$

**Case 3:** $\delta \in (0, \frac{1}{2})$. In this case, $k = j$ and $\gamma = \delta + \frac{1}{2}$, so that inequality (3.8) can be rewritten here as $q_j^{1-\delta} q_{j+1}^{\delta} \ge (\frac{1}{2} - \delta) q_j + (\frac{1}{2} + \delta) q_{j+1}$ or, equivalently, as

$$F(\delta, u) := u^{\delta} - \left(\tfrac{1}{2} - \delta\right) - \left(\tfrac{1}{2} + \delta\right) u \ge 0, \tag{3.9}$$

where

$$u := \frac{q_{j+1}}{q_j} \in [0,1]; \tag{3.10}$$

note that the conditions $x \le j_{**} + 1$ and $j \le x < j+1$ imply $j \le j_{**} + 1 \le n$ (the latter inequality takes place because, by (2.16), $j_{**} < n$); hence, $q_j \ge q_n > 0$, and thus, $u$ is correctly defined by (3.10). Moreover, because both sides of inequality (3.8) are continuous in $x$ for all $x \le n$ and hence for all $x \le j_{**} + 1$, it suffices to prove (3.8) only for $x < j_{**} + 1$, whence $j \le j_{**}$, and so, by (2.16),

$$j \le \frac{n - u_{**} q/p}{1 + u_{**} q/p};$$

the latter inequality is equivalent, in view of (2.21), to $\frac{p_{j+1}}{p_j} \ge u_{**}$, which in turn implies, in view of (2.17), that

$$\frac{q_j}{q_{j+1}} = 1 + \frac{p_j}{q_{j+1}} \le 1 + \frac{p_j}{p_{j+1}} \le 1 + \frac{1}{u_{**}} = \frac{1}{u_*},$$

whence, by (3.10), one obtains $u \ge u_*$. Therefore, the proof in Case 3, and hence the entire proof of Lemma 3.4, is now reduced to the following lemma. □

**Lemma 3.5.** *Inequality* (3.9) *holds for all* $\delta \in (0, \frac{1}{2})$ *and* $u \in [u_*, 1]$.

*Proof of Lemma 3.5.* Observe first that, for every $\delta \in (0, \frac{1}{2})$,

$$\begin{aligned} &F(\delta, 0) = -\left(\tfrac{1}{2} - \delta\right) < 0, \ F(\delta, 1) = 0, \\ &(\partial_u F)(\delta, 1) = -\tfrac{1}{2} < 0, \text{ and } F \text{ is concave in } u \in [0, 1]. \end{aligned} \tag{3.11}$$

This implies that, for every $\delta \in (0, \frac{1}{2})$, there exists a unique value

$$u(\delta) \in (0, 1)$$

such that

$$F(\delta, u(\delta)) = 0, \tag{3.12}$$

and at that

$$F(\delta, u) \geq 0 \quad \forall \delta \in \left(0, \tfrac{1}{2}\right) \ \forall u \in [u(\delta), 1], \tag{3.13}$$

and $(\partial_u F)(\delta, u(\delta))$ is strictly positive and hence nonzero. Thus, equation (3.12) defines an implicit function $(0, \frac{1}{2}) \ni \delta \mapsto u(\delta) \in (0, 1)$. Moreover, since $F$ is differentiable on $(0, \frac{1}{2}) \times (0, 1)$ and $(\partial_u F)(\delta, u(\delta)) \neq 0$ for all $\delta \in (0, \frac{1}{2})$, the implicit function theorem is applicable, so that $u(\delta)$ is differentiable in $\delta$ for all $\delta \in (0, \frac{1}{2})$. Now, differentiating both sides of equation (3.12) in $\delta$, one obtains

$$u^\delta \ln u + 1 - u + (\delta u^{\delta - 1} - \tfrac{1}{2} - \delta) u'(\delta) = 0, \tag{3.14}$$

where $u$ stands for $u(\delta)$.

Let us now show that $u(0+) = u(\frac{1}{2}-) = 0$. To that end, observe first that

$$\sup_{0 < \delta < \frac{1}{2}} u(\delta) < 1. \tag{3.15}$$

Indeed, otherwise there would exist a sequence $(\delta_n)$ in $(0, \frac{1}{2})$ such that $\varepsilon_n := 1 - u(\delta_n) \downarrow 0$. But then $u(\delta_n)^{\delta_n} = (1 - \varepsilon_n)^{\delta_n} = 1 - \varepsilon_n \delta_n + o(\varepsilon_n)$, so that (3.12) would imply

$$0 = 1 - \varepsilon_n \delta_n + o(\varepsilon_n) - \left(\tfrac{1}{2} - \delta_n\right) - \left(\tfrac{1}{2} + \delta_n\right)(1 - \varepsilon_n) = (\tfrac{1}{2} + o(1)) \varepsilon_n,$$

which would contradict the fact that $\varepsilon_n = 1 - u(\delta_n) > 0$ for all $n$.

If it were not true that $u(0+) = 0$, then there would exist a sequence $(\delta_n)$ in $(0, \frac{1}{2})$ and some $\varepsilon > 0$ such that $\delta_n \downarrow 0$ while $u(\delta_n) \to \varepsilon$. But then $u(\delta_n)^{\delta_n} \to 1$, so that equation (3.12) would imply

$$u(\delta_n) = \frac{u(\delta_n)^{\delta_n} - \left(\frac{1}{2} - \delta_n\right)}{\frac{1}{2} + \delta_n} \to 1,$$

which would contradict (3.15). Thus, $u(0+) = 0$.

Similarly, if it were not true that $u(\frac{1}{2}-) = 0$, then there would exist a sequence $(\delta_n)$ in $(0, \frac{1}{2})$ and some $\varepsilon > 0$ such that $\delta_n \uparrow \frac{1}{2}$ while $u(\delta_n) \to \varepsilon$. But then equation (3.12) would imply

$$u(\delta_n) = \frac{u(\delta_n)^{\delta_n} - \left(\frac{1}{2} - \delta_n\right)}{\frac{1}{2} + \delta_n} \to \varepsilon^{1/2},$$

which would imply $0 < \varepsilon = \varepsilon^{1/2}$, so that $\varepsilon = 1$, which would contradict (3.15). Hence, $u(\frac{1}{2}-) = 0$.

Thus, $(0,\frac{1}{2})\ni\delta\mapsto u(\delta)$ is a strictly positive continuous function, which vanishes at the endpoints 0 and $\frac{1}{2}$. Therefore, there must exist a point $\delta_*\in(0,\frac{1}{2})$ such that $u(\delta_*)\ge u(\delta)$ for all $\delta\in(0,\frac{1}{2})$. Then one must have $u'(\delta_*)=0$.

Now equation (3.14) yields

$$u(\delta_*)^{\delta_*}=-\frac{1-u(\delta_*)}{\ln u(\delta_*)},\quad\text{whence}\tag{3.16}$$

$$\delta_*=\frac{\ln\frac{1-u(\delta_*)}{-\ln u(\delta_*)}}{\ln u(\delta_*)}.\tag{3.17}$$

In the expression (3.9) for $F(\delta,u)$, replace now $u^\delta$ by the right-hand side of (3.16), and then replace $\delta$ by the right-hand side of (3.17). Then, recalling (3.12) and slightly re-arranging terms, one sees that $u=u(\delta_*)$ is a root of equation (2.15).

By Lemma 3.3, such a root of (2.15) is unique in $(0,1)$. It follows that $\max_{0<\delta<1/2}u(\delta)=u(\delta_*)=u_*=0.00505\dots$. In view of (3.13), this completes the proof of Lemma 3.5. □

*Proof of Proposition 2.7.* In this proof, we shall use Propositions 2.9 and 2.10 (which will be proved later in this paper) and the following preliminary remarks.

According to [23, Remark 13], the restriction of the tail function $Q_n$ to the set $\mathbb{Z}$ of all integers is log-concave. Therefore, the logarithm, $\ln Q_n^{\mathsf{LC}}$, of the least log-concave majorant $Q_n^{\mathsf{LC}}$ of $Q_n$ can be obtained by the linear interpolation of $\ln Q_n$ over $\mathbb{Z}$, so that

$$Q_n^{\mathsf{LC}}(x)=q_j^{1-(x-j)}q_{j+1}^{x-j}\quad\text{if}\quad j\le x\le j+1\ \&\ j\in\mathbb{Z},\tag{3.18}$$

where $q_j$ is defined by (2.21). Here and elsewhere, $0^0:=1$. Recall that the function $Q_n^{\mathsf{Lin}}$ is the linear interpolation of the function $Q_n$ over $\mathbb{Z}$, so that

$$Q_n^{\mathsf{Lin}}(x)=(1-(x-j))q_j+(x-j)q_{j+1}\quad\text{if}\quad j\le x\le j+1\ \&\ j\in\mathbb{Z}.\tag{3.19}$$

Since the function $Q_n$ is non-decreasing and left-continuous, one can note that $Q_n^{\mathsf{Lin}}\ge Q_n$ on $\mathbb{R}$; also, $Q_n^{\mathsf{Lin}}=Q_n$ on $\mathbb{Z}$.

Let us now proceed to the proof of Proposition 2.7.

**(i)** Part (i) of Proposition 2.7 follows by Lemma 3.3.

**(ii)** If $x\le j_{**}+1$ and $x$ is not in $\delta_j$ for any $j\in\mathbb{Z}\cap[j_*,n]$, then, by (2.25) and Lemma 3.4,

$$Q_n^{\mathsf{Lin,LC}}\left(x+\tfrac{1}{2}\right)=Q(x)=Q_n^{\mathsf{Lin}}\left(x+\tfrac{1}{2}\right)\le Q_n^{\mathsf{LC}}(x).\tag{3.20}$$

If $x\in\delta_j\subset(-\infty,j_{**}+1]$ for some $j\in\mathbb{Z}\cap[j_*,n]$, then, taking also into account the definition (2.24), (3.20), and the log-concavity of the function $Q_n^{\mathsf{LC}}$, one has, for $\delta$ as in (2.24),

$$\begin{aligned}Q_n^{\mathsf{Lin,LC}}\left(x+\tfrac{1}{2}\right)&=Q(x)=Q_n^{\mathsf{Interp}}(x;j)\\&=Q_n^{\mathsf{Lin}}\left(y_j+\tfrac{1}{2}\right)^{1-\delta}Q_n^{\mathsf{Lin}}\left(x_j+\tfrac{1}{2}\right)^{\delta}\le Q_n^{\mathsf{LC}}(y_j)^{1-\delta}Q_n^{\mathsf{LC}}(x_j)^{\delta}\le Q_n^{\mathsf{LC}}(x).\end{aligned}$$

Hence,

$$Q_n^{\mathsf{Lin,LC}}\left(x+\tfrac{1}{2}\right)\le Q_n^{\mathsf{LC}}(x)\tag{3.21}$$

for all $x \le j_{**}+1$ except maybe when $x \in \delta_j \cap (-\infty, j_{**}+1]$ for some $j \in \mathbb{Z} \cap [j_*, n]$ such that $\delta_j \not\subseteq (-\infty, j_{**}+1]$. The latter exceptional situation implies that $y_j < j_{**}+1 < x_j$. Hence, by Proposition 2.9,

$$j-1 < y_j < j_{**}+1 < x_j \le j+\tfrac{1}{2}, \tag{3.22}$$

whence $j < j_{**}+2$ and $j \ge j_{**}+1$, so that $j = j_{**}+1$.

It follows that (3.21) holds (at least) for all

$$x \in (-\infty, j_{**}+1] \setminus \delta_{j_{**}+1} = (-\infty, y_{j_{**}+1}] \supset (-\infty, j_{**}],$$

the latter inclusion taking place because of the first inequality in (3.22), for $j = j_{**}+1$. This completes the proof of part (ii) of Proposition 2.7.

**(iii)** Introduce

$$\tilde{Q}(x) := \begin{cases} Q_n^{\mathsf{LC}}(x) & \text{if } x \le n, \\ q_{n-1}^{n-x} q_n^{x-n+1} & \text{if } x \ge n-1. \end{cases} \tag{3.23}$$

Note that for $x \in [n-1, n]$ the expressions for $\tilde{Q}(x)$ in the two cases in (3.23) coincide with each other, which implies that the function $\tilde{Q}$ is log-concave on $\mathbb{R}$.

Now we need the following lemma, whose proof will be given a bit later.

**Lemma 3.6.** *Under condition* (2.18)*, for all* $x \in [n, n+\frac{1}{2}]$*,*

$$Q_n^{\mathsf{Lin}}\left(x+\tfrac{1}{2}\right) \le \tilde{Q}(x). \tag{3.24}$$

Under condition (2.18), it is easy to see that $j_{**}+1 \ge n$. Therefore, by Lemma 3.4, one has (3.8) and hence (3.24) for all $x \le n$.

Also, $Q_n^{\mathsf{Lin}}\left(x+\frac{1}{2}\right) = 0$ for all $x \ge n+\frac{1}{2}$.

Thus, by Lemma 3.6, $\tilde{Q}(x) \ge Q_n^{\mathsf{Lin}}\left(x+\frac{1}{2}\right)$ for all $x \in \mathbb{R}$; at that, as noted above, the function $\tilde{Q}$ is log-concave on $\mathbb{R}$. Hence, $Q_n^{\mathsf{Lin,LC}}\left(x+\frac{1}{2}\right) \le \tilde{Q}(x)$ for all $x \in \mathbb{R}$. Now part (iii) of Proposition 2.7 follows in view of (3.23), because $\tilde{Q} = Q_n^{\mathsf{LC}}$ on the interval $(-\infty, n]$. □

*Proof of Lemma 3.6.* In view of the definitions (3.19) and (3.18), one can rewrite inequality (3.24) as $q_{n-1}^{n-x} q_n^{x-n+1} \ge \left(n+\frac{1}{2}-x\right) q_n + \left(x-n+\frac{1}{2}\right) q_{n+1}$, for all $x \in [n, n+\frac{1}{2}]$. Since $q_{n+1} = 0$, the latter inequality is equivalent to

$$\ln u \le r(\alpha), \tag{3.25}$$

where $u := q_{n-1}/q_n > 1$, $\alpha := x-n \in (0, \frac{1}{2})$, and

$$r(\alpha) := \frac{\ln(\frac{1}{2}-\alpha)}{-\alpha}.$$

Since $r'(0+) = -\infty$, $r'(\frac{1}{2}-) = \infty$, and $((\ln(\frac{1}{2}-\alpha))'_\alpha)/(-\alpha'_\alpha) = 1/(\frac{1}{2}-\alpha)$ is increasing in $\alpha \in (0, \frac{1}{2})$, it follows from Proposition 3.2 that there is a unique value $\alpha_* \in (0, \frac{1}{2})$ such that the function $r$ is decreasing on $(0, \alpha_*)$ and increasing on $(\alpha_*, \frac{1}{2})$, so that $r'(\alpha_*) = 0$ and $\alpha_*$ is the point of minimum of function $r$ on $(0, \frac{1}{2})$. In fact, one has $\alpha_* = 0.3133\dots$ and $r(\alpha_*) = 5.3566\dots$ .

Therefore, inequality (3.25) can be rewritten as $u \le e^{r(\alpha_*)}$. On the other hand,

$$u = \frac{q_{n-1}}{q_n} = \frac{np^{n-1}q + p^n}{p^n} = n\,\frac{q}{p} + 1,$$

so that it suffices to check that $n\,\frac{q}{p} \le e^{r(\alpha_*)} - 1 = 211.022\dots$ . But this follows from condition (2.18). □

*Proof of Proposition 2.9.* Assume the condition $j \ge j_*$. For $j \ge n+1$, the inequalities of Proposition 2.9 immediately follow from (2.23).

It remains to consider the case $j \le n$, which implies $p_{j-1} > p_j$, by (2.22). Then it suffices to check four inequalities, $j-1 < y_j < j-\frac12 < x_j < y_{j+1}$. Indeed, the inequality $j-\frac32 < j-1$ is trivial, and the inequality $y_{j+1} \le j+\frac12$ will then follow — for $j=n$, from (2.23); and for $j=j_*,\dots,n-1$, from the inequalities $y_i < i-\frac12$ $\forall i = j_*,\dots,n$.

**(i)** Checking $j-1<y_j$. In view of definition (2.20),

$$y_j = j - \tfrac12 + \kappa_j q_j, \tag{3.26}$$

where

$$\kappa_j := \frac{1}{p_{j-1}} + \frac{\frac{1}{p_{j-1}} - \frac{1}{p_j}}{\ln\frac{p_{j-1}}{p_j}} = \frac{1 - \frac{p_{j-1}}{p_j} + \ln\frac{p_{j-1}}{p_j}}{p_{j-1}\ln\frac{p_{j-1}}{p_j}},$$

so that

$$\kappa_j < 0, \tag{3.27}$$

in view of the condition $p_{j-1} > p_j$ and the inequality $\ln u < u-1$ for $u>1$.

On the other hand, it is well known and easy to verify that the probability mass function $(p_k)$ of the binomial distribution is log-concave, so that the ratio $p_k/p_{k-1}$ is decreasing in $k$. Hence,

$$q_j = \sum_{k=j}^{n} p_k \le \sum_{k=j}^{\infty} p_j \left(\frac{p_j}{p_{j-1}}\right)^{k-j} = \frac{p_j}{1-\frac{p_j}{p_{j-1}}} =: \hat q_j.$$

Therefore, to check $j-1<y_j$, it suffices to check that $d_j := \hat y_j - (j-1) > 0$, where $\hat y_j := j - \frac12 + \kappa_j \hat q_j$ (cf. (3.26). But one can see that $d_j = f(u)/(2(u-1)\ln u)$, where $u := \frac{p_{j-1}}{p_j} > 1$ and $f(u) := 2(1-u) + (1+u)\ln u$. Thus, to check $j-1<y_j$, it suffices to show that $f(u)>0$ for $u>1$. But this follows because $f(1)=f'(1)=0$ and $f$ is strictly convex on $(1,\infty)$.

**(ii)** Checking $y_j < j-\frac12$. This follows immediately from (3.26) and (3.27).

**(iii)** Checking $j-\frac12 < x_j$. This follows because $x_j - \big(j-\frac12\big) = q_j(\ln u - 1 + 1/u)/(p_j \ln u) > 0$, where again $u := \frac{p_{j-1}}{p_j} > 1$.

**(iv)** Checking $x_j < y_{j+1}$. Let first $j \le n-1$, so that $p_j > p_{j+1} > 0$. In view of (2.20) and the obvious identity $1 + \frac{q_{j+1}}{p_j} = \frac{q_j}{p_j}$, one has

$$y_{j+1} = j - \frac12 + \frac{q_j}{p_j} + \frac{\frac{q_{j+1}}{p_j} - \frac{q_{j+1}}{p_{j+1}}}{\ln\frac{p_j}{p_{j+1}}},$$

so that the inequality $x_j < y_{j+1}$, which is being checked, can be rewritten as $q_j r_j > q_{j+1} r_{j+1}$ or, equivalently, as

$$\sum_{k=0}^{\infty} (p_{j+k} r_j - p_{j+1+k} r_{j+1}) > 0, \tag{3.28}$$

where $r_j := (\frac{1}{p_j} - \frac{1}{p_{j-1}})/\ln\frac{p_{j-1}}{p_j}$. Note that $r_j>0$, $r_{j+1}>0$, and $p_j r_j = h(v) := (1-v)/-\ln v$, where $v := p_j/p_{j-1} \in (0,1)$. By Proposition 3.1, $h(v)$ is increasing in $v\in(0,1)$. On the other hand, $v = p_j/p_{j-1}$ is decreasing in $j$, by the mentioned log-concavity of $(p_j)$. It follows that $p_j r_j$ is decreasing in $j$. Because of this and the same log-concavity, $\frac{p_{j+k}r_j}{p_{j+1+k}r_{j+1}} > \frac{p_j r_j}{p_{j+1}r_{j+1}} > 1 \quad \forall k = 0, 1, \ldots$, which yields (3.28).

Finally, in the case when $j = n \geq j_*$, the inequality $x_j < y_{j+1}$ follows from (2.19) and (2.23), because then

$$x_n = n + \frac{1}{2} + \frac{\dfrac{p_n}{p_{n-1}} - 1}{\ln \dfrac{p_{n-1}}{p_n}} < n + \frac{1}{2} = y_{n+1}.$$

□

*Proof of Proposition 2.10.*

**Step 1.** Here we observe that the function $Q$ defined in (2.25) is continuous. Indeed, the function $x \mapsto Q_n^{\mathsf{Lin}}\left(x + \frac{1}{2}\right)$ is defined and continuous everywhere on $\mathbb{R}$. On the other hand, for every integer $j \geq j_*$, the function $x \mapsto Q_n^{\mathsf{Interp}}(x; j)$ is defined and continuous on the interval $\delta_j$; moreover, it continuously interpolates on the interval $\delta_j$ between the values of the function $x \mapsto Q_n^{\mathsf{Lin}}\left(x + \frac{1}{2}\right)$ at the endpoints, $y_j$ and $x_j$, of the interval $\delta_j$. Also, the intervals $\delta_j$ with $j \in \mathbb{Z} \cap [j_*, n]$ are pairwise disjoint. Thus, the function $Q$ is continuous everywhere on $\mathbb{R}$.

**Step 2.** Here we show that the function $Q$ is log-concave. To that end, introduce

$$\ell_j(x) := \ln\left((\tfrac{1}{2} + j - x)q_j + (\tfrac{1}{2} - j + x)q_{j+1}\right) \quad \forall j \in \mathbb{Z}, \quad \text{so that}$$

$$(3.29) \qquad (j - \tfrac{1}{2} \leq x < j + \tfrac{1}{2} \ \&\ j \leq n) \implies \ln Q_n^{\mathsf{Lin}}\left(x + \tfrac{1}{2}\right) = \ell_j(x);$$

here, the condition $j \leq n$ provides for both sides of the equality in (3.29) to be defined. One can check (which is better done using Mathematica or similar software) the basic relations

$$(3.30) \qquad j \in \mathbb{Z} \cap [j_*, n] \implies \ell'_{j-1}(y_j) = \ell'_j(x_j) = \frac{\ell_j(x_j) - \ell_{j-1}(y_j)}{x_j - y_j};$$

these relations do not rely on the fact that the $q_j$'s pertain to a binomial distribution, but only on general relations: $p_{j-1} > p_j > 0$, $p_i = q_i - q_{i+1}\ \forall i$, and $q_{j+1} \geq 0$, as well as the inequalities $(j-1) - \frac{1}{2} < y_j < (j-1) + \frac{1}{2}$ and $j - \frac{1}{2} < x_j \leq j + \frac{1}{2}$, which follow by Proposition 2.9 and ensure that $\ell_j$ and $\ell_{j-1}$ are defined and differentiable in neighborhoods of $x_j$ and $y_j$, respectively. Using the latter relations together with (2.24) and (3.29), one has

$$(3.31) \qquad \frac{\mathrm{d}}{\mathrm{d}x} \ln Q_n^{\mathsf{Interp}}(x; j) = \frac{\ell_j(x_j) - \ell_{j-1}(y_j)}{x_j - y_j} \quad \forall x \in \delta_j\ \forall j \in \mathbb{Z} \cap [j_*, n].$$

Moreover, for all integer $j \leq n$,

$$\ell'_j\left(j - \tfrac{1}{2}\right) = \frac{q_{j+1} - q_j}{q_j} = \frac{-p_j}{q_j} \quad \text{and} \quad \ell'_{j-1}\left(j - \tfrac{1}{2}\right) = \frac{-p_{j-1}}{q_j}.$$

Hence and by (2.22), for every integer $j \leq n$ one has

$$(3.32) \qquad \ell'_j\left(j - \tfrac{1}{2}\right) \leq \ell'_{j-1}\left(j - \tfrac{1}{2}\right) \iff j \leq j_* - 1.$$

In view of (3.29), the function $x \mapsto \ln Q_n^{\mathsf{Lin}}\left(x + \frac{1}{2}\right)$ is concave on the interval $[j - \frac{1}{2}, j + \frac{1}{2}]$ for every integer $j \leq n$. Note also that $\bigcup_{j \leq j_* - 1}\left[j - \frac{1}{2}, j + \frac{1}{2}\right] = \left(-\infty, j_* - \frac{1}{2}\right]$. Hence, by (3.32) and (3.29), the function

$$(3.33) \qquad x \mapsto \ln Q_n^{\mathsf{Lin}}\left(x + \tfrac{1}{2}\right) \quad \text{is concave on the interval } \left(-\infty, j_* - \tfrac{1}{2}\right].$$

In addition to the open intervals $\delta_j = (y_j, x_j)$, introduce the closed intervals

$$\Delta_j := [x_j, y_{j+1}]$$

for integer $j \geq j_*$. Then, by Proposition 2.9 and (2.23), the intervals $\Delta_j$ are each nonempty,

$$\delta_{j_*} \cup \Delta_{j_*} \cup \delta_{j_*+1} \cup \Delta_{j_*+1} \cup \cdots \cup \delta_n \cup \Delta_n = \big(y_{j_*}, n + \tfrac{1}{2}\big], \tag{3.34}$$

and

$$\delta_{j_*} < \Delta_{j_*} < \delta_{j_*+1} < \Delta_{j_*+1} < \cdots < \delta_n < \Delta_n.$$

Thus, the intervals $\delta_j$ and $\Delta_j$ with $j \in \mathbb{Z} \cap [j_*, n]$ form a partition of the interval $\big(y_{j_*}, n + \frac{1}{2}\big]$. Moreover, for every $j \in \mathbb{Z} \cap [j_*, n]$, by Proposition 2.9, $\Delta_j \subseteq [j - \frac{1}{2}, j + \frac{1}{2}]$, and so, by (3.29), the function

$$x \mapsto \ln Q_n^{\mathsf{Lin}}\big(x + \tfrac{1}{2}\big) \quad \text{is concave on the interval } \Delta_j. \tag{3.35}$$

Also, by (2.24), for every $j \in \mathbb{Z} \cap [j_*, n]$, the function

$$x \mapsto \ln Q_n^{\mathsf{Interp}}(x; j) \text{ is concave (in fact, affine) on the interval } \delta_j. \tag{3.36}$$

By the definition of $Q$ in (2.25), for all $x \in \mathbb{R}$ and all $j \in \mathbb{Z} \cap [j_*, n]$, one has

$$Q(x) = \begin{cases} Q_n^{\mathsf{Lin}}\big(x + \frac{1}{2}\big) & \text{if } \ x \leq y_{j_*}, \\ Q_n^{\mathsf{Interp}}(x; j) & \text{if } \ x \in \delta_j \ \& \ j \in \mathbb{Z} \cap [j_*, n], \\ Q_n^{\mathsf{Lin}}\big(x + \frac{1}{2}\big) & \text{if } \ x \in \Delta_j \ \& \ j \in \mathbb{Z} \cap [j_*, n], \\ 0 = Q_n^{\mathsf{Lin}}\big(x + \frac{1}{2}\big) & \text{if } \ x \geq n + \frac{1}{2}. \end{cases} \tag{3.37}$$

Note also that, by Proposition 2.9, $y_{j_*} \leq j_* - \frac{1}{2}$, so that $(-\infty, y_{j_*}] \subseteq (-\infty, j_* - \frac{1}{2}]$. Now it follows from (3.37), (3.33), (3.35), and (3.36) that the function $\ln Q$ is concave on each of the disjoint adjacent intervals

$$(-\infty, y_{j_*}], \delta_{j_*}, \Delta_{j_*}, \delta_{j_*+1}, \Delta_{j_*+1}, \ldots, \delta_n, \Delta_n, \tag{3.38}$$

whose union is the interval $(-\infty, n + \frac{1}{2}]$. Moreover, it follows from the continuity of $Q$ (established in Step 1) and formulas (3.30), (3.29), and (3.31) that the function $\ln Q$ is differentiable at all the endpoints $y_{j_*}, x_{j_*}, y_{j_*+1}, x_{j_*+1}, \ldots, y_n, x_n$ of the intervals (3.38) except the right endpoint $y_{n+1} = n + \frac{1}{2}$ of the interval $\Delta_n$.

Therefore, the function $\ln Q$ is concave on the interval $(-\infty, n + \frac{1}{2})$. On the other hand, $\ln Q = -\infty$ on the interval $[n + \frac{1}{2}, \infty)$. Thus, it is proved that the function $\ln Q$ is concave everywhere on $\mathbb{R}$.

**Step 3.** Here we show that

$$Q(x) \geq Q_n^{\mathsf{Lin}}\big(x + \tfrac{1}{2}\big) \tag{3.39}$$

for all real $x$. In view of (3.37) and (3.34), it suffices to check (3.39) for $x \in \delta_j$ with $j \in \mathbb{Z} \cap [j_*, n]$. By Proposition 2.9, $\delta_j \subseteq (j - 1, j + \frac{1}{2}] \subset [j - \frac{3}{2}, j + \frac{1}{2}]$, for every $j \in \mathbb{Z} \cap [j_*, n]$.

By (3.29), the function $x \mapsto \ln Q_n^{\mathsf{Lin}}\big(x + \frac{1}{2}\big) = \ell_j(x)$ is concave on the interval $[j - \frac{1}{2}, j + \frac{1}{2}]$, for every integer $j \leq n$. Hence, (3.30) and (2.24) imply that, for all $x \in \delta_j \cap [j - \frac{1}{2}, j + \frac{1}{2}]$ with $j \in \mathbb{Z} \cap [j_*, n]$,

$$\begin{aligned} \ln Q_n^{\mathsf{Lin}}\big(x + \tfrac{1}{2}\big) &= \ell_j(x) \leq \ell_j(x_j) + \ell_j'(x_j)(x - x_j) \\ &= \frac{x_j - x}{x_j - y_j}\, \ell_{j-1}(y_j) + \frac{x - y_j}{x_j - y_j}\, \ell_j(x_j) = \ln Q_n^{\mathsf{Interp}}(x; j) = \ln Q(x), \end{aligned}$$

so that one has (3.39) for all $x \in \delta_j \cap [j - \frac{1}{2}, j + \frac{1}{2}]$ with $j \in \mathbb{Z} \cap [j_*, n]$. Similarly (using inequality $\ell_{j-1}(x) \le \ell_{j-1}(y_j) + \ell'_{j-1}(y_j)(x - y_j)$) it can be shown that (3.39) takes place for all $x \in \delta_j \cap [j - \frac{3}{2}, j - \frac{1}{2}]$ with $j \in \mathbb{Z} \cap [j_*, n]$. This completes Step 3.

**Step 4.** Here we show that, if $\tilde{Q}$ is a log-concave function on $\mathbb{R}$ such that

$$\tilde{Q}(x) \ge Q_n^{\mathsf{Lin}}\left(x + \tfrac{1}{2}\right) \quad \forall x \in \mathbb{R}, \tag{3.40}$$

then $\tilde{Q} \ge Q$ on $\mathbb{R}$. In view of (3.37), it suffices to check that $\tilde{Q} \ge Q$ on $\delta_j$ for every $j \in \mathbb{Z} \cap [j_*, n]$. But, by (3.40), one has $\tilde{Q}(y_j) \ge Q_n^{\mathsf{Lin}}\left(y_j + \frac{1}{2}\right)$ and $\tilde{Q}(x_j) \ge Q_n^{\mathsf{Lin}}\left(x_j + \frac{1}{2}\right)$. Hence, taking into account the log-concavity of $\tilde{Q}$ and (2.24) and (3.37), one has, for all $x \in \delta_j$ with $j \in \mathbb{Z} \cap [j_*, n]$ and $\delta$ as in (2.24),

$$\tilde{Q}(x) \ge \tilde{Q}\left(y_j\right)^{1-\delta} \tilde{Q}\left(x_j\right)^{\delta} \ge Q_n^{\mathsf{Lin}}\left(y_j + \tfrac{1}{2}\right)^{1-\delta} Q_n^{\mathsf{Lin}}\left(x_j + \tfrac{1}{2}\right)^{\delta} = Q_n^{\mathsf{Interp}}(x; j) = Q(x).$$

The facts established in Steps 2, 3, and 4 imply that the function $Q$ is indeed the least log-concave majorant of the function $x \mapsto Q_n^{\mathsf{Lin}}\left(x + \frac{1}{2}\right)$. Thus, Proposition 2.10 is proved. □